\documentclass[11pt,a4paper]{article}
\usepackage[T1]{fontenc}
\usepackage[utf8]{inputenc}
\usepackage{lmodern}
\usepackage[margin=25mm]{geometry}
\usepackage{amsmath,amssymb,amsthm,booktabs,microtype}
\usepackage{graphicx,tikz}
\usepackage[numbers,sort&compress]{natbib}
\usepackage{etoolbox}
\usepackage[section]{placeins}
\usepackage[colorlinks=true,linkcolor=blue,citecolor=blue,urlcolor=blue]{hyperref}
\AtBeginEnvironment{thebibliography}{\raggedright}
\newtheorem{theorem}{Theorem}
\newtheorem{lemma}[theorem]{Lemma}
\newtheorem{corollary}[theorem]{Corollary}
\newtheorem{proposition}[theorem]{Proposition}
\theoremstyle{remark}

\newcommand{\R}{\mathbb R}

\newcommand{\Col}{\operatorname{Col}}
\newcommand{\Span}{\operatorname{span}}
\newcommand{\calI}{\mathcal I}

\newcommand{\supp}{\operatorname{supp}}
\hypersetup{
 pdftitle={Repairing actuator schedules under graph constraints},
 pdfauthor={Dawid E. Kedzierski and Piotr Krason},
 pdfsubject={Preprint: actuator scheduling under graph constraints},
 pdfkeywords={actuator scheduling, basis exchanges, compatible bases,
 controllability, graph constraints, reachability profiles}}
\hypersetup{
  colorlinks=true,
  allcolors=black
}
\title{Repairing actuator schedules under graph constraints}
\author{Dawid E. K{\k{e}}dzierski \and Piotr Kraso{\'n}}
\date{}

\begin{document}
\maketitle
\begin{NoHyper}

\begin{abstract}
We consider a discrete-time linear system in which each actuator is assigned
to one time slot. A graph $H$ records which pairs of actuators cannot share
a slot. We seek to remove conflicts by exchanging actuators between slots,
while preserving all tail reachability spaces. These spaces describe the
terminal corrections available from each time slot onward. We give a
sufficient condition based on the number of exchanges that preserve the
required bases. The condition guarantees repair from every compatible
starting schedule. Two bounds for checking it use failed basis tests and
dependent column selections. Each accepted exchange reduces the conflict
count, so the procedure uses at most as many exchanges as there are initial
conflicts. An explicit family with unequal rank increments admits repair
when $3N>4\Delta(H)+4$, where $\Delta(H)$ is the largest number of forbidden
partners of any actuator. We also give a family of two-slot examples with
only one forbidden partner per actuator and two initial conflicts. The basis
constraints force every strict improvement to move arbitrarily many
actuators, although one exchange suffices without those constraints.
A second example has a connected conflict graph and one initial conflict.
\end{abstract}

\section{Introduction}

Actuator scheduling determines when each actuator in a control system is
used. We consider schedules in which each actuator is assigned to one time
slot. Some pairs of actuators are not allowed to share a slot. We record these forbidden pairs as the
edges of a simple graph $H$ on the actuator labels. For each actuator,
count the other actuators that cannot share its slot. We denote the largest
of these counts by $\Delta(H)$ (the maximum degree of $H$).
We start from a schedule that preserves the system's tail reachability
spaces. The question is whether exchanges between slots can remove all
conflicts while preserving those spaces.

K{\k{e}}dzierski and Kraso\'n \cite{KK2026} describe one-use schedules
with this reachability profile by compatible bases. They also give an exact
search procedure. We assume that a starting schedule has already been
found. Our task is to remove its graph conflicts.

Actuator scheduling has also been studied under other constraints.
Ballotta, Joseph and Rahul Thete \cite{BJT2024} construct pointwise-sparse
schedules with a controllability guarantee. Chamon, Amice and Ribeiro
\cite{CAR2022} study quadratic control cost under operational constraints
expressed by intersecting independence systems. General families of
admissible input supports are considered in \cite{JNKV2021}.
Here we fix the horizon, assign each actuator to one slot, and preserve
every tail reachability space during repair.

Without the basis constraints, we need to divide the actuators into groups
of prescribed sizes, with no forbidden pair in a group. This is a
graph-packing problem. Let $r_k$ be the size of slot $k$,
$\rho=\max_k r_k$, and $d$ the total number of actuators. The
Sauer--Spencer theorem \cite{SS1978} gives the sufficient condition
$d>2\Delta(H)(\rho-1)$. List packing \cite{GKMY2016} also allows forbidden
vertex assignments. Our basis constraints apply to the entire selection
of actuators in a slot. The proof uses the graph-packing exchange count,
but excludes exchanges that fail a basis test. Related work includes
coloring under graph and matroid constraints \cite{AMPZ2026} and coloring
with prescribed class sizes \cite{Birken2026}. In our problem, the
independence test depends on the slot, and every slot size is fixed.

Our main result gives a sufficient condition for repair from any compatible
starting schedule. It requires enough basis-preserving exchange partners
for each actuator. One bound counts exchanges ruled out by the basis
constraints; another counts dependent column selections. We give a family
with unequal rank increments for which the first bound can be less restrictive.
We then construct a conflict graph whose edges have no shared endpoints
(a matching). The basis constraints can force every strict improvement to
move arbitrarily many actuators, although one exchange fixes the graph
conflicts when those constraints are removed. A second construction gives
a connected graph with one initial conflict.

\section{Schedules and the reachability flag}

Consider
\begin{equation}\label{eq:system}
 x_k=Ax_{k-1}+Gu_k,\qquad k=1,\ldots,N,
\end{equation}
where $A\in\R^{d\times d}$,
$G=[g_1,\ldots,g_d]\in\R^{d\times d}$, and $u_k\in\R^d$.
The actuator labels form the set $E=\{1,\ldots,d\}$.
For an input vector $u\in\R^d$, we write
\[
 \supp(u)=\{i\in E:u_i\ne0\}
\]
for the set of indices of its nonzero components.
A one-use schedule is a map $\tau:E\to\{1,\ldots,N\}$.
Write $I_k=\tau^{-1}(k)$ and require $\supp(u_k)\subseteq I_k$.
Thus only actuators assigned to slot $k$ may have nonzero inputs in that
slot. The assignment is made once per actuator; its amplitude may be zero.
The terminal direction associated with actuator $i$ is
$A^{N-\tau(i)}g_i$.

For $k=1,\ldots,N$, let
\begin{equation}\label{eq:flag}
 W_k=\Col[\,A^{N-k}G\mid\cdots\mid AG \mid G\,],\qquad W_{N+1}=0.
\end{equation}
Assume $W_1=\R^d$, and set
\begin{gather}
 V_k=W_k/W_{k+1},\qquad r_k=\dim V_k,\notag\\
 v_{k,i}=A^{N-k}g_i+W_{k+1}.\label{eq:quotients}
\end{gather}
Thus $\sum_k r_k=d$. A partition
$\calI=(I_1,\ldots,I_N)$ of $E$ is \emph{compatible} if
\begin{equation}\label{eq:compatible}
 \{v_{k,i}:i\in I_k\}\text{ is a basis of }V_k
 \quad\text{for every }k.
\end{equation}
In particular, $|I_k|=r_k$. Let $\cal P$ denote the family of all
compatible partitions of $E$. We assume $\cal P\ne\varnothing$.

Condition \eqref{eq:compatible} is equivalent to requiring that the terminal
directions assigned to slots $k,\ldots,N$ form a basis of $W_k$ for every
$k$. Indeed, descending induction on $k$ combines a basis of $W_{k+1}$
with representatives of a basis of $W_k/W_{k+1}$. The converse follows by
passing to the quotient. This is the flag interpretation of compatible
schedules developed in \cite{KK2026}. In particular, the selected terminal
matrix is invertible, so any terminal state can be reached from any initial
state by solving for the assigned input amplitudes.

For any one-use schedule $\tau$, define the space of terminal corrections
available from slot $k$ onward by
\begin{equation}\label{eq:scheduled-tail}
 \mathcal R_k(\tau)=
 \Span\{A^{N-\tau(i)}g_i:\tau(i)\geq k\}.
\end{equation}
\paragraph{A fixed compatible schedule}
Fix a compatible schedule $\tau$ and a time slot $k$.
Keep the initial state and all inputs in slots $1,\ldots,k-1$ fixed.
Let $b_k$ be the terminal state obtained if all later inputs are zero.

First, keep the slot assignments of $\tau$ fixed and vary only the
remaining input amplitudes allowed by this schedule. The resulting
terminal states form the set
\[
 b_k+\mathcal R_k(\tau)=b_k+W_k,
\]
where the equality follows from compatibility.

For comparison, keep the same initial state and past inputs, but allow
every actuator to be used in every remaining slot $k,\ldots,N$.
The resulting terminal states also form $b_k+W_k$.
Thus the fixed compatible schedule preserves all terminal corrections
available from slot $k$ onward.

\paragraph{Terminal controllability without compatibility}
Consider a separate example with $A=G=I_2$ and $N=2$.
Here $W_1=W_2=\R^2$ and $(r_1,r_2)=(0,2)$.
A compatible schedule therefore assigns both actuators to slot~2.

Instead, choose the schedule
\[
 \widehat\tau(1)=1,\qquad \widehat\tau(2)=2.
\]
This schedule is not compatible. Its terminal matrix is nevertheless
\[
 [Ag_1\mid g_2]=I_2,
\]
so it is terminally controllable when both input amplitudes are free.
Once the slot-1 amplitude has been fixed, however, only actuator~2
remains available. It can change only the second coordinate, and hence
\[
 \mathcal R_2(\widehat\tau)
 =\Span\{g_2\}
 \subsetneq W_2=\R^2.
\]
This schedule therefore gives full-horizon controllability but does
not preserve the last-slot reachability space.
In the repair results below, we always start from a compatible schedule.

Let $H=(E,F)$ be a simple conflict graph, with vertex set $E$ (the
actuator labels) and edge set $F$ (the forbidden pairs).
Thus $\{i,j\}\in F$ means that actuators $i$ and $j$ cannot share a slot.
For an actuator $v$, the number of other actuators forbidden to share its slot is
$\deg_H(v)=|\{u\in E:\{u,v\}\in F\}|$. We use
\[
 \Delta=\Delta(H)=\max_{v\in E}\deg_H(v)
\]
for the largest of these counts. In graph terminology, $\deg_H(v)$ is the
degree of $v$, and $\Delta(H)$ is the maximum degree of $H$. Define
\begin{equation}\label{eq:potential}
 \kappa(\calI)=\sum_{k=1}^N e_H(I_k),
\end{equation}
where $e_H(S)$ is the number of edges with both ends in $S$.
A schedule satisfies the graph constraints precisely when $\kappa=0$.
Set $\rho=\max_k r_k$. The case $\Delta=0$ requires no repair.
All subsequent bounds involving $\Delta$ concern $\Delta\geq1$.

\section{Exchange repair}

The exchanges below change the planned slot assignments before the
schedule is executed. Fix a compatible partition
$\calI=(I_1,\ldots,I_N)$.
For each slot $k\in\{1,\ldots,N\}$, the vectors
$\{v_{k,i}:i\in I_k\}$ form the current basis of $V_k$.
For every actuator label $j\in E=\{1,\ldots,d\}$, the vector
$v_{k,j}=A^{N-k}g_j+W_{k+1}$, where $g_j$ is column $j$ of $G$,
belongs to $V_k$ and has a unique expansion in this basis:
\begin{equation}\label{eq:coeff}
 v_{k,j}=\sum_{i\in I_k}c^{(k)}_{ij}v_{k,i}.
\end{equation}
In this equation, $k$ and $j$ are fixed, and only $i$ is summed.
The index $k$ identifies the slot and its quotient space $V_k$.
In $v_{k,j}$, the second index $j$ is the label of the actuator whose
direction is being expanded. This label can belong to any slot of
$\calI$. The vector $v_{k,j}$ describes its proposed use in slot $k$,
regardless of its current assignment $\tau(j)$.
The summation index $i$ runs through the actuator labels in $I_k$,
which index the current basis vectors of $V_k$.

Each $c^{(k)}_{ij}\in\R$ is a scalar coordinate.
Its first lower index $i\in I_k$ specifies the basis vector $v_{k,i}$
that it multiplies. Its second lower index $j\in E$ specifies the vector
$v_{k,j}$ being expanded. The parenthesized superscript $(k)$ identifies
the slot in which these coordinates are calculated.
All these coordinates refer to the basis selected by the fixed
partition $\calI$. More explicitly, one could write
$c^{(k)}_{ij}(\calI)$; we suppress this dependence to keep the notation
short. The coefficients may change when the partition changes.

For a fixed $k$, the coefficients form an $r_k\times d$ matrix
$C^{(k)}=(c^{(k)}_{ij})$, where $r_k=|I_k|=\dim V_k$.
The rows are indexed by the labels $i\in I_k$, and the columns by
$j=1,\ldots,d$ in this order. Choosing an order on $I_k$ fixes the
row order.
For example, if $I_k=\{2,5\}$ and $j=7$, then
\[
 v_{k,7}=c^{(k)}_{2,7}v_{k,2}+c^{(k)}_{5,7}v_{k,5}.
\]
Here $i$ takes the two values $2$ and $5$, while $j$ remains $7$.
With row order $(2,5)$, the coefficients $c^{(k)}_{2,7}$ and
$c^{(k)}_{5,7}$ occupy the first and second rows of column $7$.
The row labels are actuator labels, whose numerical values need not
equal their row positions.
If $j\in I_k$, the expanded vector is itself a basis vector:
$c^{(k)}_{ij}=1$ for $i=j$ and $c^{(k)}_{ij}=0$ for the other
$i\in I_k$. If $I_k$ is empty, then $V_k=\{0\}$ and both sides
of~\eqref{eq:coeff} are zero, with an empty sum on the right.

Suppose $a\in I_k$ and $b\in I_\ell$, where $k\ne\ell$.
We exchange their slot assignments: actuator $a$ moves to slot $\ell$,
and actuator $b$ moves to slot $k$. The resulting partition, denoted by
$\calI^{a\leftrightarrow b}$, has slots
$I_k'=(I_k\setminus\{a\})\cup\{b\}$ and
$I_\ell'=(I_\ell\setminus\{b\})\cup\{a\}$.
The primes denote the slots after the exchange. For every other slot
index $p\in\{1,\ldots,N\}\setminus\{k,\ell\}$, we set $I_p'=I_p$.
Thus every slot keeps its size.
Although $b$ is initially assigned to slot $\ell$, its proposed use in
slot $k$ is tested using $v_{k,b}$ in $V_k$.
\begin{proposition}[Compatible exchanges]\label{prop:exchange}
For any $\calI\in\cal P$, $a\in I_k$ and $b\in I_\ell$ with $k\ne\ell$,
\begin{equation}\label{eq:exchange}
 \mathcal{I}^{a\leftrightarrow b}\in\mathcal{P}
 \iff c^{(k)}_{a\, b} c^{(\ell)}_{b\,a} \neq 0
\end{equation}
\end{proposition}

\begin{proof}
In~\eqref{eq:exchange}, $c^{(k)}_{a\,b}$ is read from row $a$, column $b$
of $C^{(k)}$, using the original basis indexed by $I_k$.
The coefficient $c^{(\ell)}_{b\,a}$ is read from row $b$, column $a$
of $C^{(\ell)}$, using the original basis indexed by $I_\ell$.
Thus the two coefficients come from the coordinate matrices for two
different slots. Both are computed before the exchange.
We prove the equivalence by checking the two changed bases separately.
In slot $k$, the exchange removes $v_{k,a}$ and inserts $v_{k,b}$.
Separating the coefficient of the removed vector in~\eqref{eq:coeff}
gives
\[
 v_{k,b}=c^{(k)}_{a\,b}v_{k,a}
 +\sum_{i\in I_k\setminus\{a\}}c^{(k)}_{ i\, b}v_{k,i}.
\]
Here $a$ and $b$ are the two fixed exchange labels. 
The vectors indexed by $I_k\setminus\{a\}$ remain in the slot.
If $c^{(k)}_{a\,b}=0$, the inserted vector lies in their span.
The updated selection is then dependent and cannot be a basis of $V_k$.
This proves that $c^{(k)}_{a\,b}\ne0$ is necessary.

Conversely, if $c^{(k)}_{a\,b}\ne0$, we can recover the removed vector:
\[
 v_{k,a}=\frac{1}{c^{(k)}_{a\,b}}\left(
 v_{k,b}-\sum_{i\in I_k\setminus\{a\}}c^{(k)}_{i\, b}v_{k,i}
 \right).
\]
Thus the updated vectors span every vector of the original basis,
and hence all of $V_k$. There are still exactly $r_k=\dim V_k$ of
them, so they form a basis. This proves sufficiency as well.

In slot $\ell$, the removed label is $b$ and the inserted label is $a$.
The relevant expansion in its original basis is
\[
 v_{\ell,a}=c^{(\ell)}_{b\,a}v_{\ell,b}
 +\sum_{i\in I_\ell\setminus\{b\}}c^{(\ell)}_{i\, a}v_{\ell,i}.
\]
The same two arguments show that this updated selection is a basis
of $V_\ell$ if and only if $c^{(\ell)}_{ba}\ne0$.
The order of the indices is reversed because slot $\ell$ removes $b$
and receives $a$.

If the exchanged partition is compatible, both updated selections are
bases, so both coefficients must be nonzero. Conversely, if both
coefficients are nonzero, both updated selections are bases; all other
quotient bases are unchanged, so the entire partition is compatible. This gives exactly~\eqref{eq:exchange}.
\end{proof}

Recall that $\cal P$ is the family of partitions satisfying
the basis condition~\eqref{eq:compatible}.
Fix $\calI\in\cal P$ and choose an actuator $i\in E$.
Let $\tau$ be the schedule associated with this partition, and put
$k=\tau(i)$. Thus $i\in I_k$ and $I_{\tau(i)}=I_k$.
We collect all labels that can be exchanged with $i$ while preserving
compatibility:
\begin{equation}\label{eq:partners}
 S_{\calI}(i)=\{j\in E\setminus I_{\tau(i)}:
              \calI^{i\leftrightarrow j}\in\cal P\}.
\end{equation}
The subscript $\calI$ specifies the entire current partition, and the
argument $i$ specifies the actuator whose exchange partners are listed.
Both are fixed in this definition; $j$ runs through the candidate labels.
The set $E\setminus I_{\tau(i)}$ consists of the actuators assigned to
slots other than that of $i$. It excludes $i$ itself and every other
actuator in its current slot.
We retain $j$ exactly when both updated sets of quotient vectors
are linearly independent.

To test a candidate $j\in E\setminus I_k$, set $\ell=\tau(j)$.
Then $\ell\ne k$, and~\eqref{eq:exchange} gives
\[
 j\in S_{\calI}(i)
 \iff c^{(k)}_{ij}c^{(\ell)}_{ji}\ne0.
\]
Thus $j$ is retained exactly when both basis replacements pass their
tests. Here $k$ remains fixed, whereas $\ell=\tau(j)$ can vary with the
candidate $j$. Every candidate is tested using the same original
partition $\calI$ and its coordinate matrices.

The elements of $S_{\calI}(i)$ are actuator labels, and
$|S_{\calI}(i)|$ is the number of available compatible exchange partners.
For example, if $I_k=\{2,5\}$, $i=2$, and the tests succeed exactly for
the outside labels $3$ and $7$, then
\[
 S_{\calI}(2)=\{3,7\},\qquad |S_{\calI}(2)|=2.
\]
These are two alternative single exchanges from the same partition.
After one is performed, the partition changes and its exchange-partner
sets may change as well. If no candidate passes both tests, then
$S_{\calI}(i)=\varnothing$ and $|S_{\calI}(i)|=0$.

Membership in $S_{\calI}(i)$ concerns only compatibility with the bases.
The conflict graph $H$ does not enter this test, and the conflict count
$\kappa$ may increase, decrease, or remain unchanged after such an
exchange. The repair procedure will require, in addition, that the
chosen exchange strictly decrease $\kappa$.

We need a bound that holds for every compatible partition and every
choice of actuator. Set
\begin{equation}\label{eq:mu}
 \mu=\min_{\calI\in\cal P}\ \min_{i\in E}|S_{\calI}(i)|.
\end{equation}
For a fixed partition, the inner minimum is the smallest number of
exchange partners available to any actuator. The outer minimum takes
the least of these numbers over all compatible partitions.
Thus every actuator has at least $\mu$ compatible exchange partners,
regardless of the compatible partition currently used. The value of
$\mu$ may be zero.

The matrices $A,G$ and the horizon $N$ determine the quotient vectors,
the compatible partitions, and the sets $S_{\calI}(i)$.
No graph condition enters these definitions. Consequently, changing
the conflict graph $H$ while keeping $A,G,N$ fixed does not change $\mu$.
The graph determines which of the compatible exchanges reduce conflicts.

Recall that $r_k=|I_k|$ for a compatible partition and that
$\rho=\max_{1\leq p\leq N}r_p$ is the largest slot size.
Also, $\Delta=\Delta(H)$ is the maximum degree of the conflict graph:
each actuator has at most $\Delta$ forbidden partners.
The next proof compares the uniform lower bound $\mu$ on compatible
partners with the number of candidates that the graph can exclude.

\begin{theorem}[Exchange repair]\label{thm:repair}
If
\begin{equation}\label{eq:mucondition}
 \mu>(2\Delta-1)\rho-2\Delta,
\end{equation}
then every compatible partition with a conflict admits a compatible
exchange that strictly decreases $\kappa$. Starting from any
$\calI^{(0)}\in\cal P$, repeated improving exchanges therefore produce a
compatible partition satisfying the graph constraints in at most
$\kappa(\calI^{(0)})$ exchanges.
\end{theorem}

\begin{proof}
Fix $\calI\in\cal P$ with $\kappa(\calI)>0$.
Choose a conflict $\{v,u\}\in F$ with $v,u\in I_k$.
We seek a partner $w\in S_{\calI}(v)$ whose exchange with $v$ removes
this conflict and creates no new one.
For $S\subseteq E$, let $N_H(S)$ be the set of labels adjacent to at
least one label in $S$. Define
\begin{gather*}
 X=\bigcup_{z\in N_H(\{v\})}
       \bigl(I_{\tau(z)}\setminus\{z\}\bigr),\\
 Y=N_H(I_k\setminus\{v\}),\qquad T=X\cup Y\cup\{u\}.
\end{gather*}
Consider a candidate partner outside $I_k$.
For each neighbor $z$ of $v$, the set $X$ contains all other labels in
the slot of $z$. If the candidate is one of these other labels, the
exchange puts $v$ in a slot with its neighbor $z$.
If the candidate is $z$ itself, then $z$ leaves that slot; this is why
$z$ is omitted from that term.
The set $Y$ excludes candidates adjacent to a member of
$I_k\setminus\{v\}$, since these members remain in slot $k$.

We first bound these two sets:
\begin{align*}
 |X|&\leq
 \sum_{z\in N_H(\{v\})}(r_{\tau(z)}-1)
 \leq\Delta(\rho-1),\\
 |Y|&\leq
 \sum_{z\in I_k\setminus\{v\}}\deg_H(z)
 \leq\Delta(r_k-1).
\end{align*}
The first inequality in each line bounds a union by the sum of the
sizes of its sets. For $X$, there are at most $\Delta$ neighbors $z$,
and each contributes at most $\rho-1$ labels. For $Y$, there are
$r_k-1$ retained labels, each with at most $\Delta$ neighbors.

The chosen conflict gives an overlap that is important for the count.
Since $u$ is a neighbor of $v$ in slot $k$, the term $z=u$ in $X$
contains $I_k\setminus\{u\}$. In particular, $v\in X$.
Also, $u\in I_k\setminus\{v\}$ and $v$ is its neighbor, so $v\in Y$.
Thus $v\in X\cap Y$, and adding $u$ ensures $I_k\subseteq T$.
The set $\{u\}$ adds at most one label, while $X$ and $Y$ have at
least the label $v$ in common. Therefore
\[
 |T|\leq |X|+|Y|\leq \Delta(\rho+r_k-2).
\]

Every label in $I_k$ is already in $T$. Hence the number of excluded
candidates outside the source slot is
$|T\setminus I_k|=|T|-r_k$. Subtracting $r_k$ from the preceding bound
and collecting terms gives
\begin{align*}
 |T\setminus I_k|
 &\leq\Delta\rho+(\Delta-1)r_k-2\Delta\\
 &\leq(2\Delta-1)\rho-2\Delta.
\end{align*}
The second inequality uses $r_k\leq\rho$ and the nonnegative coefficient
$\Delta-1$. Here $\Delta\geq1$ because the graph contains the chosen
conflict edge.

By definition, $|S_{\calI}(v)|\geq\mu$, and every member of
$S_{\calI}(v)$ lies outside $I_k$.
Assumption~\eqref{eq:mucondition} therefore gives
$|S_{\calI}(v)|>|T\setminus I_k|$.
Consequently there is a partner $w\in S_{\calI}(v)\setminus T$.
Membership in $S_{\calI}(v)$ ensures that the exchanged partition
remains compatible.

Put $\ell=\tau(w)$; then $\ell\ne k$. If $v$ had a neighbor
$z\in I_\ell\setminus\{w\}$, then
$w\in I_{\tau(z)}\setminus\{z\}\subseteq X$, contrary to $w\notin T$.
Thus inserting $v$ in slot $\ell$ creates no conflict.
Also, $w\notin Y$ ensures that inserting $w$ in slot $k$ creates no
conflict with $I_k\setminus\{v\}$.
The pair $\{v,w\}$ lies in different slots before and after the exchange,
and all pairs not involving either label keep their assignments.
The old conflict $\{v,u\}$ is removed, so
\[
 \kappa(\calI^{v\leftrightarrow w})\leq\kappa(\calI)-1.
\]
The same argument applies to every compatible partition with positive
conflict count, so the procedure cannot stop while a conflict remains.
Since $\kappa$ is a nonnegative integer and decreases by at least one
at each accepted exchange, it reaches zero in at most
$\kappa(\calI^{(0)})$ exchanges.
\end{proof}

\paragraph{Repair procedure}
Given a compatible starting partition, repeat the following operations
while $\kappa>0$: choose an endpoint $v$ of an internal graph edge; scan
actuators in other slots; accept the first exchange passing
\eqref{eq:exchange} for which $\kappa$ decreases. Under
Theorem~\ref{thm:repair}, this scan always succeeds, including when the
conflicted endpoint is chosen arbitrarily.

There are at most $d\,\kappa(\calI^{(0)})$ candidate tests. Each candidate
needs two quotient-basis tests and a graph-conflict update involving the
two exchanged vertices. Given quotient matrices, these tests use ordinary
linear algebra. Neither $\mu$ nor a complete enumeration of compatible
partitions is needed by the repair procedure itself. Obtaining or
certifying a suitable lower bound for $\mu$ is a separate task.

\subsection{A bound on algebraically blocked partners}

For a compatible partition and $v\in I_k$, let
\begin{equation}\label{eq:blocked}
 a_{\calI}(v)=d-r_k-|S_{\calI}(v)|.
\end{equation}
This counts partners in other slots for which at least one of the two
quotient-basis tests fails. A bound valid for every actuator and compatible
partition gives the following sufficient condition.

\begin{corollary}[Bounded algebraic exclusions]\label{cor:blocked}
Suppose $a_{\calI}(v)\leq a$ for every $\calI\in\cal P$ and $v\in E$.
If
\begin{equation}\label{eq:blockedcondition}
 d>2\Delta(\rho-1)+a,
\end{equation}
then the repair procedure terminates without conflicts from every
compatible starting partition in at most $\kappa(\calI^{(0)})$ exchanges,
preserving every tail reachability space after each exchange.
\end{corollary}

\begin{proof}
Equation~\eqref{eq:blocked} gives
$|S_{\calI}(v)|\geq d-r_k-a\geq d-\rho-a$.
Thus \eqref{eq:blockedcondition} implies
$\mu\geq d-\rho-a>(2\Delta-1)\rho-2\Delta$.
Apply Theorem~\ref{thm:repair}.
\end{proof}

The number $a$ bounds failed exchange partners for one actuator in a given
partition. It can be smaller than the total number of dependent column
selections. A bound valid for all compatible partitions must be proved
separately. Section~4 gives an example with $a=1$.

\subsection{A condition in terms of dependent selections}

For each quotient $V_k$, let $b_k$ count the sets $T\subseteq E$ with
$|T|=r_k$ for which the vectors $\{v_{k,i}:i\in T\}$ are linearly
dependent. Define
\begin{equation}\label{eq:bad}
 b=\sum_{k=1}^N b_k.
\end{equation}
The empty selection for $r_k=0$ is independent and contributes zero.
Selections with the same labels in different quotients are counted
separately.

\begin{lemma}\label{lem:bad}
For every $\calI\in\cal P$ and $v\in I_k$,
\[
 a_{\calI}(v)\leq b,\qquad |S_{\calI}(v)|\geq d-r_k-b.
\]
In particular, $\mu\geq d-\rho-b$.
\end{lemma}

\begin{proof}
There are $d-r_k$ possible partners outside $I_k$. A failed exchange
produces a dependent selection in either the departing or the receiving
slot. Fix a quotient index $\ell$ and a selection $T$ counted by $b_\ell$.
If $\ell=k$,
the equation $T=(I_k\setminus\{v\})\cup\{w\}$ determines at most one
partner $w$. If $\ell\ne k$, the equation
$T=(I_\ell\setminus\{w\})\cup\{v\}$ again determines at most one
partner. Thus each pair $(\ell,T)$ excludes at most one partner. A union
bound excludes at most $b$ partners in total.
\end{proof}

\begin{corollary}[Rank-count certificate]\label{cor:bad}
If
\begin{equation}\label{eq:badcondition}
 d>2\Delta(\rho-1)+b,
\end{equation}
then the repair procedure terminates without conflicts from every
compatible starting partition, preserving every tail reachability space
after each exchange.
\end{corollary}

\begin{proof}
Lemma~\ref{lem:bad} gives the uniform bound $a=b$.
Apply Corollary~\ref{cor:blocked}.
\end{proof}

Lemma~\ref{lem:bad} shows that $b$ bounds the true maximum number of
failed partners over all compatible partitions and actuators.
Whenever the condition using $b$ holds, the condition using that true
maximum also holds. The role of $b$ is that it can be obtained from
the individual quotient matrices without listing all compatible
partitions. This is useful when the dependent column selections are
known explicitly. A smaller uniform bound $a$ requires a separate
proof valid for every compatible partition and every actuator.

For $b=0$, the degree condition is exactly the Sauer--Spencer bound
\cite{SS1978}; the additional term accounts for dependent selections.
Listing all subsets to compute $b$ can still be expensive.
The condition is sufficient: a failed check does not rule out a
conflict-free compatible schedule.

\section{Unequal rank increments}

We compare two estimates of the number of exchange partners ruled out
by the basis conditions. For this family, counting failed exchanges
directly gives the uniform bound $a=1$, while counting all dependent
selections gives $b=N-2$. The quotient dimensions are unequal, with
maximum $\rho=3$.

Let $N\geq3$ and $d=3N-3$. Take three Jordan chains with respective
lengths $N,N-1,N-2$, with bases
\[
 a_0,\ldots,a_{N-1};\quad
 b_0,\ldots,b_{N-2};\quad
 z_0,\ldots,z_{N-3}.
\]
The operator $A$ sends each vector to its successor and the last vector
of each chain to zero. Define
\begin{equation}\label{eq:familyG}
 g_i=\alpha_i a_0+\beta_i b_0+\gamma_i z_0.
\end{equation}
For each actuator $i=1,\ldots,d$, let
$c_i:=(\alpha_i,\beta_i,\gamma_i)^{\mathsf T}\in\R^3$.
The index $i$ labels the actuator. The three entries of $c_i$ are the
coefficients of $g_i$ at the starting vectors $a_0,b_0,z_0$, in that
order. We choose these coefficient vectors as follows:
\begin{gather}
 c_1=(1,1,0)^{\mathsf T},\quad c_2=(2,3,0)^{\mathsf T},\quad
 c_3=(3,4,0)^{\mathsf T},\notag\\
 c_i=(1,i,i^2)^{\mathsf T}\quad(i=4,\ldots,d).\label{eq:Cfamily}
\end{gather}
Thus the first three columns of the input matrix $G=[g_1,\ldots,g_d]$
are
\[
 g_1=a_0+b_0,\quad g_2=2a_0+3b_0,\quad g_3=3a_0+4b_0.
\]
This choice gives $c_3=c_1+c_2$ and hence $g_3=g_1+g_2$.
The first three columns are therefore linearly dependent, while every
pair among them is linearly independent.

The first quotient has coordinate columns $\alpha_i$, the second has
columns $(\alpha_i,\beta_i)^{\mathsf T}$, and all remaining quotients
have columns $c_i$. Therefore
\begin{equation}\label{eq:familyr}
 (r_1,\ldots,r_N)=(1,2,3,\ldots,3).
\end{equation}

\begin{proposition}[Explicit family]\label{prop:family}
For the family \eqref{eq:familyG}--\eqref{eq:Cfamily},
$\rho=3$, $b=N-2$ and $a_{\calI}(v)\leq1$ for every compatible partition
and actuator. In particular, $\mu\geq d-4$.
For every simple conflict graph satisfying
$3N>4\Delta(H)+4$, any compatible starting partition can be repaired by
the procedure of Theorem~\ref{thm:repair}.
\end{proposition}

\begin{proof}
Every $\alpha_i$ is nonzero. The ratios $\beta_i/\alpha_i$ are
$1,\allowbreak3/2,\allowbreak4/3,\allowbreak4,\allowbreak5,\ldots,d$, so every pair of columns in the second
quotient is independent. The only dependent triple of full coefficient
columns is $\{1,2,3\}$, for which $c_3=c_1+c_2$.

To verify the latter claim, two of the first three columns together with
any $c_i$, $i\geq4$, have nonzero determinant. Three columns with labels
at least four form a Vandermonde matrix. Finally, for
$c=(\alpha,\beta,0)^{\mathsf T}$ among the first three columns and
$4\leq s<t$, the determinant with $c_s,c_t$ is
\[
 (t-s)\bigl(\alpha st-\beta(s+t)\bigr)\ne0,
\]
because $\beta/\alpha\leq3/2$ and $st/(s+t)\geq2$.
Thus $b_1=b_2=0$ and $b_k=1$ for $k\geq3$.

A compatible partition exists explicitly: take $I_1=\{1\}$,
$I_2=\{2,3\}$ and
$I_k=\{3k-5,3k-4,3k-3\}$ for $k\geq3$. The quotient ranks give
\eqref{eq:familyr} and full reachability.

Fix any compatible schedule and choose one actuator $v$.
Consider exchanging $v$ with each actuator assigned to another slot.
We show that at most one of these exchanges violates the basis condition.
Write $\calI$ for the partition of this schedule, with $v\in I_k$.

The determinant calculations above show that the only dependent triple
in each three-dimensional quotient is formed by the columns of actuators
$1,2,3$. Let $Q=\{1,2,3\}\subseteq E$ be the set of these three actuators.
A slot of size three therefore fails the basis condition exactly when
its actuator set is $Q$. All selections of the required size for slots
of size one or two are bases, as also shown above.

Consider a candidate $w\in I_\ell$ with $\ell\ne k$.
Exchanging $v$ and $w$ changes the two slots to
$(I_k\setminus\{v\})\cup\{w\}$ and
$(I_\ell\setminus\{w\})\cup\{v\}$.
All other slots keep their bases. Hence this exchange is incompatible
exactly when one of these two updated label sets equals $Q$.
Here we count only failures of the basis condition over all candidates
$w$ outside $I_k$.

If $v\in Q$, the updated slot $k$ cannot equal $Q$, because it no longer
contains $v$. Only slot $\ell$ can become $Q$. For this to happen, its
two retained labels must be exactly $Q\setminus\{v\}$.
These two fixed labels can lie together in at most one slot of the
partition. If that slot has size three, the only failing partner $w$
is its third label, which is removed when $v$ enters. Thus there is
at most one failing partner in this case.

If $v\notin Q$, the updated slot $\ell$ contains $v$, so it cannot
equal $Q$. Only slot $k$ can become $Q$. This requires $|I_k|=3$ and
the two retained labels in $I_k\setminus\{v\}$ to belong to $Q$.
Then $w$ must be the unique third label of $Q$ missing from that pair.
Since $I_k$ and $v$ are fixed, this again excludes at most one partner.

Both cases give $a_{\calI}(v)\leq1$. The schedule and actuator were
arbitrary, so this bound holds for every compatible schedule and every
actuator. It is independent of $N$.
Among the $d-r_k$ labels outside $I_k$, at most one fails the basis tests, so
$|S_{\calI}(v)|\geq d-r_k-1$.
Since $r_k\leq3$, taking the minima over $\calI$ and $v$ gives
$\mu\geq d-4=3N-7$.
We may therefore use $a=1$ and $\rho=3$ in
Corollary~\ref{cor:blocked}. Its condition becomes $d>4\Delta+1$,
or equivalently $3N>4\Delta+4$, as claimed.
If $\Delta=0$, there are no graph conflicts and no repair is needed.
\end{proof}

The count $b=N-2$ includes the same triple $Q$ once for each of the
$N-2$ three-dimensional quotients. Using this count,
condition~\eqref{eq:badcondition} reads
$3N-3>4\Delta+(N-2)$, equivalently $N\geq2\Delta+1$.
The direct exchange count instead uses $a=1$, regardless of how many
quotients contain this dependent triple.
Because $a=1\leq b=N-2$, the condition using $b$ always implies the
condition using $a=1$. The calculation in the next paragraph shows that
the condition using $a=1$ can hold while the condition using $b$ fails.

For a concrete comparison, take $N=6$ and a graph with maximum degree
$\Delta=3$. Then $d=15$, $\rho=3$ and $b=4$.
In condition~\eqref{eq:blockedcondition}, the right-hand side is
$2\Delta(\rho-1)+a=13$. Since $15>13$, repair is guaranteed from every
compatible starting partition.
In condition~\eqref{eq:badcondition}, the right-hand side is
$2\Delta(\rho-1)+b=16$. Since $15<16$, this condition is not satisfied
and gives no guarantee for the same system and graph.

\subsection{Exact-arithmetic verification}

For $N=3$, the system has dimension six and increments $(1,2,3)$.
Exactly one last-slot triple is dependent. Hence there are
$19\cdot3=57$ compatible partitions. All $76$ matchings on six labels
were enumerated, including the empty matching; these are all simple
graphs of maximum degree at most one. Repair succeeded for all
$57\cdot76=4332$ pairs of starting partition and graph. The minimum
number of compatible exchange partners was $\mu=2$, and at most two
exchanges were needed.
For $F=\{12,34,56\}$, one run is shown in
Table~\ref{tab:small}; the exchanges are $1\leftrightarrow3$ and then
$1\leftrightarrow5$.

\begin{table}[htbp]
\centering\small
\begin{tabular}{ccccc}
\toprule
Iteration & $I_1$ & $I_2$ & $I_3$ & $\kappa$\\
\midrule
0 & $\{3\}$ & $\{1,2\}$ & $\{4,5,6\}$ & 2\\
1 & $\{1\}$ & $\{2,3\}$ & $\{4,5,6\}$ & 1\\
2 & $\{5\}$ & $\{2,3\}$ & $\{1,4,6\}$ & 0\\
\bottomrule
\end{tabular}
\vspace{4pt}
\caption{Repair of a six-dimensional schedule with conflict-graph edges $\{12,34,56\}$. Each $I_k$ lists the actuator labels assigned to slot $k$; $\kappa$ counts conflict edges whose ends share a slot.}
\label{tab:small}
\end{table}

All 4920 compatible partitions at $N=4$ were also enumerated. Every
vertex had at most one algebraically blocked partner, and $\mu=5=d-4$.
The graph exclusion count in the repair proof was checked for every graph
and set partition on five labels.

Larger instances used $N=5,6,7$ with degree caps $2,3,4$,
respectively. For each cap, 100 graphs were generated by randomly
ordering all possible edges and accepting an edge if both current degrees
were below the cap. Each graph received an independently sampled
compatible starting partition. All 300 runs ended without conflicts.
The largest exchange counts were $4,7,8$, respectively. These are
sampled checks, not exhaustive counts. These checks used integer
determinants and rational rank calculations, with seed 20261023.

\section{Limits of strict descent}

A compatible schedule need not satisfy the graph constraints. For example,
if every pair of actuators is forbidden, each nonempty slot can contain
only one actuator. Even when a conflict-free compatible schedule exists,
small updates need not reduce the conflict count.
Proposition~\ref{prop:matching} shows how the basis conditions can force
an improving update to move many actuators.

\begin{proposition}[Matching conflict graph]\label{prop:matching}
For every integer $q\geq2$, there is a two-slot system with $4q$
actuators, a matching conflict graph of maximum degree one, and a
compatible partition with two conflicts such that every compatible
reassignment strictly reducing the conflict count moves at least $2q$
actuators. A compatible reassignment moving exactly $2q$ actuators
removes all conflicts. If the basis constraints are omitted, a single
exchange of two actuators removes both conflicts while preserving the
slot sizes.
\end{proposition}

Figure~\ref{fig:matching} illustrates the initial partition for $q=2$.

\begin{figure}[htbp]
\centering
\begin{tikzpicture}[
  x=1cm,y=1cm,
  actuator/.style={circle,draw=black,fill=white,inner sep=1pt,
    minimum size=6mm,font=\small},
  pair/.style={densely dashed,gray!75,line width=0.55pt},
  edge/.style={black,line width=0.65pt},
  conflict/.style={red!70!black,line width=1.25pt}
]
\filldraw[fill=gray!7,draw=gray!40,rounded corners=3pt]
  (0,0.45) rectangle (5.6,2.1);
\filldraw[fill=gray!7,draw=gray!40,rounded corners=3pt]
  (0,-1.8) rectangle (5.6,-0.45);
\node[font=\small] at (-0.45,1) {$I_1$};
\node[font=\small] at (-0.45,-1) {$I_2$};
\node[actuator] (x1) at (0.7,1) {$x_1$};
\node[actuator] (x2) at (2.1,1) {$x_2$};
\node[actuator] (y3) at (3.5,1) {$y_3$};
\node[actuator] (y4) at (4.9,1) {$y_4$};
\node[actuator] (y1) at (0.7,-1) {$y_1$};
\node[actuator] (y2) at (2.1,-1) {$y_2$};
\node[actuator] (x3) at (3.5,-1) {$x_3$};
\node[actuator] (x4) at (4.9,-1) {$x_4$};
\draw[pair] (x1) -- (y1);
\draw[pair] (x2) -- (y2);
\draw[pair] (x3) -- (y3);
\draw[pair] (x4) -- (y4);
\draw[edge] (y1) -- (x2);
\draw[edge] (y3) -- (x4);
\draw[conflict] (y4) .. controls (4.9,2.25) and (0.7,2.25) .. (x1);
\draw[conflict] (y2) .. controls (2.1,-1.8) and (3.5,-1.8) .. (x3);
\end{tikzpicture}
\caption{The matching construction for $q=2$, with
$s=(0,0,1,1)$. Solid edges form the conflict graph $H$; the two bold
edges are conflicts inside the slots. Dashed links join actuators with
identical columns; each slot must contain one actuator from each such pair.}
\label{fig:matching}
\end{figure}
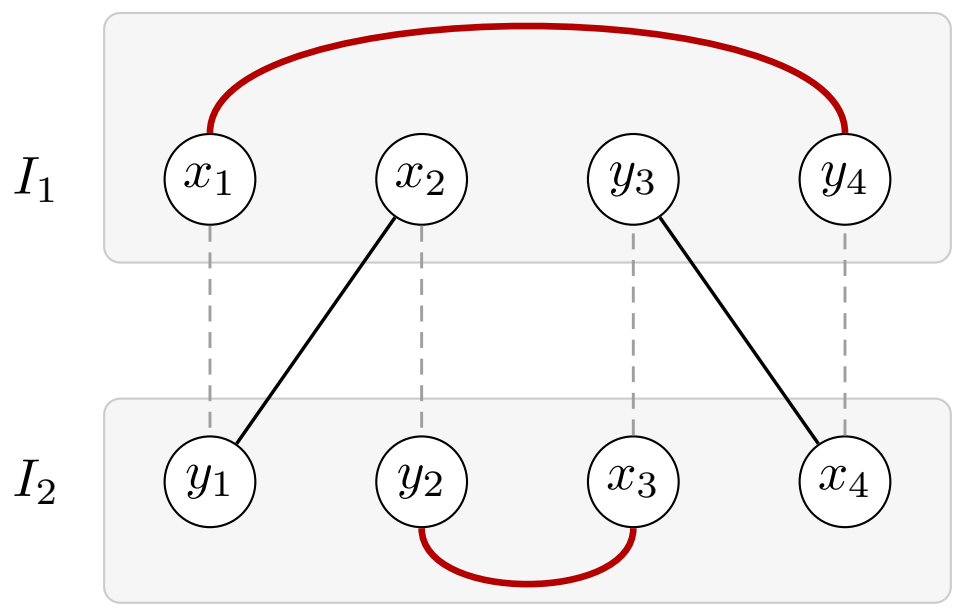

\begin{proof}
Take $N=2$, put $m=2q$, and label the $2m$ actuators
$x_i,y_i$, $1\leq i\leq m$.
With all $x_i$ columns preceding the $y_i$ columns, set
\begin{gather*}
 A=\begin{pmatrix}0&I_m\\0&0\end{pmatrix},\qquad
 G=\begin{pmatrix}0\\C\end{pmatrix},\qquad C=[\,I_m\mid I_m\,].
\end{gather*}
The columns of $G$ span the bottom $m$ coordinates, and those of $AG$
span the top $m$ coordinates. Thus $W_1=\R^{2m}$,
$W_2=\Col G$, and $r_1=r_2=m$.
In either quotient, the columns for $x_i$ and $y_i$ are both the $i$th
standard basis vector. A selection of $m$ columns is therefore a basis
exactly when it contains one actuator from every pair $(x_i,y_i)$.
Hence a compatible partition places one member of each pair in $I_1$
and the other in $I_2$.

Record these choices by $s_i=0$ when $x_i\in I_1$ and $s_i=1$ when
$x_i\in I_2$. The actuator $y_i$ is always in the other slot.

Let
\[
 F=\{\{y_i,x_{i+1}\}:1\leq i\leq m\},\qquad x_{m+1}=x_1.
\]
Each actuator belongs to exactly one edge, so the graph is a matching.
The edge $\{y_i,x_{i+1}\}$ is a conflict exactly when its ends share a
slot. Since $y_i$ is opposite $x_i$, this happens exactly when
$s_i\ne s_{i+1}$. Consequently,
\begin{equation}\label{eq:matching-conflicts}
 \kappa(s)=\sum_{i=1}^{m}\mathbf 1_{\{s_i\ne s_{i+1}\}},
 \qquad s_{m+1}=s_1.
\end{equation}
The indicator in this sum is one when its condition holds and zero
otherwise. Initially choose $s_1=\cdots=s_q=0$ and
$s_{q+1}=\cdots=s_m=1$. The only changes of value occur at $i=q$
and $i=m$, where the cyclic sequence returns to $s_1$.
Thus the two conflicts are $\{y_q,x_{q+1}\}$ in $I_2$ and
$\{y_m,x_1\}$ in $I_1$.

In a cyclic binary sequence, every change from zero to one is matched
by a change from one to zero. Hence the sum in
\eqref{eq:matching-conflicts} is even. A strict decrease from two
conflicts must therefore give zero conflicts, which requires all bits
to be zero or all bits to be one. Either target changes exactly $q$
bits of the initial sequence. Changing one bit moves both actuators
of its pair to the opposite slots. Thus a compatible strict
improvement moves $2q$ actuators, and either target attains this bound.

If we impose only the graph constraints, we can instead exchange
$y_q$ and $y_m$. Both conflict edges then have ends in different
slots. No other edge is affected, since the graph is a matching.
Each slot loses one actuator and gains one, so its size stays $m$.
However, $I_1$ now contains both $x_q$ and $y_q$, while $I_2$ contains
both $x_m$ and $y_m$. Each of these pairs has identical columns in its
quotient. Both updated selections are therefore linearly dependent,
so this exchange does not preserve compatibility.
\end{proof}

A sequence of small compatible exchanges can still repair this
schedule if the conflict count is allowed to stay unchanged at some
steps. Exchange $x_q$ with $y_q$, then $x_{q-1}$ with $y_{q-1}$,
and continue down to $x_1$ with $y_1$. Each step changes one of the
initial zero bits to one and keeps one member of every pair in each
slot. After each of the first $q-1$ steps, the cyclic bit sequence still
has one nonempty block of zeros and one of ones, so $\kappa=2$.
The last step makes all bits one and gives $\kappa=0$.
The obstruction in Proposition~\ref{prop:matching} therefore concerns
one update required to decrease $\kappa$ strictly.

A complementary construction retains a connected conflict graph and
only one initial conflict.

\begin{samepage}
\begin{proposition}[Connected conflict graph]\label{prop:barrier}
For every integer $q\geq2$, there is a two-slot system with $4q$
actuators and a connected bipartite conflict graph of maximum degree
$q+1$, together with a compatible partition having one conflict, such
that every compatible reassignment strictly reducing conflicts moves at
least $2q$ actuators. A reassignment moving exactly $2q$ actuators
removes all conflicts.
\end{proposition}
\end{samepage}

\begin{proof}
Take $N=2$ and partition the actuator labels into $L_1,R_1,L_2,R_2$,
each of size $q$. For $j=1,2$, $L_j$ and $R_j$ denote the left and right
sets of actuators in the $j$th complete bipartite subgraph.
Order the columns by these four groups and set
\begin{gather*}
 A=\begin{pmatrix}0&I_{2q}\\0&0\end{pmatrix},\qquad
 G=\begin{pmatrix}0\\C\end{pmatrix},\\
 C=\begin{pmatrix}I_q&I_q&0&0\\0&0&I_q&I_q\end{pmatrix}.
\end{gather*}
As before, the columns of $G$ and $AG$ span complementary coordinate subspaces, so
$r_1=r_2=2q$. For $j=1,2$, pair the $i$th actuator
of $L_j$ with the $i$th actuator of $R_j$, for $1\leq i\leq q$.
The members of each pair have identical columns. As in the preceding
proof, compatibility requires exactly one member of each pair in each
slot.

Put an edge between every actuator of $L_1$ and every actuator of $R_1$,
and do the same for $L_2$ and $R_2$. Add one bridge from a fixed
actuator of $L_1$ to a fixed actuator of $R_2$.
The graph is connected and has maximum degree $q+1$.
Here and below, the two components are those left after deleting the bridge.

Initially let $I_1=L_1\cup R_2$ and $I_2=R_1\cup L_2$.
Every edge within either bipartite component joins different slots.
Both ends of the bridge are in $I_1$, so there is exactly one conflict.

A compatible reassignment exchanges the members of some pairs.
Let $s$ be the number of such pairs in the first component and $t$
the number in the second, with $0\leq s,t\leq q$.
In the first component, slot $I_1$ then has $q-s$ actuators from $L_1$
and $s$ from $R_1$. Since every pair between these groups is an edge,
this gives $s(q-s)$ conflicts in $I_1$ and the same number in $I_2$.
The second component similarly gives $2t(q-t)$ conflicts. Thus
\begin{equation}\label{eq:barrier}
 \kappa=2s(q-s)+2t(q-t)+\varepsilon,
\end{equation}
where $\varepsilon=1$ if the bridge ends share a slot and
$\varepsilon=0$ otherwise. A strict decrease from one conflict must
give zero conflicts. If $0<s<q$, both $s$ and $q-s$ are positive,
so the first term is positive. The same argument applies to $t$.
Thus $s,t\in\{0,q\}$:
each component must keep all its original assignments or reverse all
of them. The bridge ends lie in different slots exactly when one
component is reversed and the other is unchanged.
This exchanges all $q$ pairs of one component, moving $2q$ actuators.
It preserves the bases and removes every conflict, so the lower
bound is attained.
\end{proof}

In this connected construction, the graph alone forces the same number
of moved actuators. A conflict-free assignment must put all of $L_j$
in one slot and all of $R_j$ in the other, for each component.
To separate the bridge ends, one whole component must reverse its
initial assignments, moving $2q$ actuators.
Proposition~\ref{prop:matching} shows a different effect: there the graph
alone is repaired by one exchange, and the basis conditions force the
larger update.

Both propositions count the actuators whose slot assignments change
in one update that strictly decreases $\kappa$. They allow sequences
of small updates in which $\kappa$ stays unchanged or increases
temporarily; the sequence following Proposition~\ref{prop:matching}
is one such repair.
Exhaustive enumeration of all $2^{2q}$ compatible partitions for
$q=2,\ldots,8$ confirmed \eqref{eq:matching-conflicts},
\eqref{eq:barrier}, and both minimum numbers of moved actuators.

For any proposed reassignment preserving slot sizes, let
$D_k\subseteq I_k$ be the labels leaving slot $k$ and
$P_k\subseteq E\setminus I_k$ the labels entering it.
Then $|D_k|=|P_k|$. Use the coordinate matrix $C^{(k)}$ of the current
partition from~\eqref{eq:coeff}. The square submatrix
$C^{(k)}[D_k,P_k]$ has rows labeled by $D_k$ and columns labeled by
$P_k$.
Pass to the quotient of $V_k$ by the span of the retained vectors,
whose labels are $I_k\setminus D_k$. This quotient has dimension
$|D_k|$, and the inserted vectors have coordinate matrix
$C^{(k)}[D_k,P_k]$.

Hence the updated slot has a basis exactly when
\[
 \det C^{(k)}[D_k,P_k]\ne0.
\]
Checking this condition in every changed slot tests whether the
whole reassignment is compatible. It also applies to larger updates
when no single exchange gives a compatible strict improvement.

\section{Discussion}

Under the stated conditions, repair removes all graph conflicts from
any compatible starting schedule in at most $\kappa(\calI^{(0)})$
accepted exchanges. Each exchange preserves every tail reachability
space, and hence every terminal correction available after earlier
input amplitudes have been fixed.

The family in Section~4 shows the benefit of counting failed partners
directly: $a=1$ stays fixed while $b=N-2$ grows with the horizon.
In the matching construction, the basis conditions force a large
update whenever the conflict count must decrease strictly.
Allowing exchanges that keep the count unchanged permits a sequence
of small compatible updates. In the connected construction, the graph
alone forces the large update.
Further work could combine the rank conditions with bounds on control
energy and matrix conditioning.

\texttt{TailBasisGraphRepair}~\cite{TailBasisGraphRepair2026} is a supplementary
Wolfram Language package.
It takes a compatible schedule and a conflict graph as input and extends
the workflow of TailBasis~\cite{TailBasis2026}
by displaying compatible exchanges, successive schedules, and the conflict
count after each step. The package is distributed separately.

\end{NoHyper}


\bigskip
\begingroup
\small
\noindent\textbf{Dawid E. K{\k{e}}dzierski}\\
Institute of Mathematics, University of Szczecin\\
ul. Wielkopolska 15, 70-451 Szczecin, Poland\\
\nolinkurl{dawidedmundkedzierski@gmail.com}

\medskip
\noindent\textbf{Piotr Kraso{\'n}}\\
Department of Mathematics, Ariel University, Ariel 40700, Israel\\
\nolinkurl{piotrkras26@gmail.com}
\endgroup

\end{document}